\documentclass[a4paper,11pt]{amsart}
\usepackage{graphicx}
\usepackage{mathptmx}
\usepackage{amsmath}
\usepackage{amssymb}
\usepackage{enumitem}
\usepackage{xcolor}
\usepackage{xparse}
\NewDocumentCommand{\eulerian}{omm}
 {%
  \genfrac<>{0pt}{}{#2}{#3}%
  \IfValueT{#1}{_{\!#1}}%
 }

 \newmuskip\pFqmuskip

\newcommand*\pFq[6][8]{%
  \begingroup 
  \pFqmuskip=#1mu\relax
  \mathchardef\normalcomma=\mathcode`,
  \mathcode`\,=\string"8000
  \begingroup\lccode`\~=`\,
  \lowercase{\endgroup\let~}\pFqcomma
  {}_{#2}F_{#3}{\left(\genfrac..{0pt}{}{#4}{#5}\bigg|#6\right)}%
  \endgroup
}
\newcommand{\pFqcomma}{{\normalcomma}\mskip\pFqmuskip}

\newtheorem{theorem}{Theorem}

\newtheorem{corollary}[theorem]{Corollary}

\newtheorem{remark}[theorem]{Remark}

\begin{document}

\title[Some results on degenerate harmonic numbers and degenerate Fubini polynomials]{Some results on degenerate harmonic numbers and degenerate Fubini polynomials}

\author{Taekyun  Kim}
\address{Department of Mathematics, Kwangwoon University, Seoul 139-701, Republic of Korea}
\email{tkkim@kw.ac.kr}

\author{DAE SAN KIM}
\address{Department of Mathematics, Sogang University, Seoul 121-742, Republic of Korea}
\email{dskim@sogang.ac.kr}

\subjclass[2010]{11B73; 11B83}
\keywords{degenerate harmonic numbers; degenerate hyperharmonic numbers; degenerate Fubini polynomials; degenerate $r$-Fubini polynomials}

\begin{abstract}
In recent years, some degenerate versions of quite a few special numbers and polynomials are introduced and investigated by means of various methods. The aim of this paper is to study some results on degenerate harmonic numbers, degenerate hyperharmonic numbers, degenerate Fubi polynomials and degenerate $r$-Fubini polynomials from a general identity which is valid for any two formal power series and involves the degenerate $r$-Stirling numbers of the second kind.
\end{abstract}

\maketitle

\section{Introduction}
It was Carlitz who initiated the study of degenerate versions of some special numbers and polynomials, namely the degenerate Bernoulli and degenerate Euler numbers (see [4]). This pioneering work regained interests of some mathematicians, various versions of many special numbers and polynomials were investigated and some interesting arithmetical and combinatorial results were found along with some of their applications (see [9-11,13-22] and the references therein).  \par
The aim of this paper is to study some results on degenerate harmonic numbers, degenerate hyperharmonic numbers, degenerate Fubi polynomials and degenerate $r$-Fubini polynomials from a general identity which is valid for any two formal power series and involves the degenerate $r$-Stirling numbers of the second kind. \par
The outline of this paper is as follows. In Section 1, we recall the degenerate exponentials, their compositional inverses called the degenerate logarithms, the degenerate Stirling numbers of the first kind and the degenerate Stirling numbers of the second kind. We mention the degenerate $r$-Stirling numbers of the first kind, the degenerate Stirling $r$-Stirling numbers of the second kind and the degenerate unsigned $r$-Stirling numbers of the first kind as further generalizations of aforementioned numbers. We recall the reader of the degenerate Bell polynomials, the degenerate $r$-Bell polynomials, the degenerate harmonic numbers, the degenerate hyperharmonic numbers and the degenerate Fubini polynomials. Section 2 is the main results of this paper. In Theorem 1, we derive an identity of operators involving the degenerate $r$-Stirling numbers of the second kind. A general identity valid for any two formal power series is derived in Theorem 2 by using Theorem 1, where again the degenerate $r$-Stirling numbers of the second kind appear. Applying the result in Theorem 2 with $f(x)=x^m$ and $g(x)=\frac{1}{1-x}$, we get an identity on the degenerate $r$-Fubini polynomials which are slight generalization of the degenerate Fubini polynomials (see \eqref{36}). In Theorem 4, we deduce a differential equation satisfied by the degenerate Fubini polynomials with the help of a recurrence relation for the degenerate Stirling numbers of the second kind. In Theorem 5, the degenerate hyperharmonic numbers are expressed in terms of the degenerate unsigned $r$-Stirling numbers of the first kind. Theorem 6 shows that the $k$-th derivative of a certain formal power series is equal to an expression involving the degenerate harmonic number $H_{k,\lambda}$ and the degenerate logarithm whose value at 0 is simply $k!H_{k,\lambda}$. In Corollary 7, we derive a relation between the degenerate hyperharmonic numbers and the degenerate harmonic numbers from Theorem 6 and some other previous results. Finally, by applying Theorem 2 with $f(x)=x^m$ and $g(x)=-\frac{1}{1-x}\log_{\lambda}(1-x)$ we obtain an identity expressing a power series with coefficients given by the product of the degenerate harmonic numbers and the degenerate falling factorials in terms of the degenerate $r$-Stirlings of the second, the degenerate harmonic numbers and the degenerate logarithms. Finally, we conclude this paper in Section 3.

\par
For any nonzero $\lambda\in\mathbb{R}$, the degenerate exponentials are defined by 
\begin{equation}
e_{\lambda}^{x}(t)=\sum_{n=0}^{\infty}\frac{(x)_{n,\lambda}}{n!}t^{n},\quad e_{\lambda}^{1}(t)=e_{\lambda}(t),\quad (\mathrm{see}\ [12,13,14]),\label{1}
\end{equation}
where the degenerate falling factorials are given by
\begin{equation}
(x)_{0,\lambda}=1,\quad (x)_{n,\lambda}=x(x-\lambda)\cdots (x-(n-1)\lambda),\quad (n\ge 1).\label{2}
\end{equation}
Note that $\lim_{\lambda\rightarrow 0}e_{\lambda}^{x}(t)=\sum_{n=0}^{\infty}\frac{x^{n}}{n!}t^{n}=e^{xt}$. \par 
Let $\log_{\lambda}t$ be the compositional inverse of $e_{\lambda}(t)$, called the degenerate logarithm, satisfying $e_{\lambda}(\log_{\lambda}t)=\log_{\lambda}(e_{\lambda}(t))=t$. \par 
Then we have 
\begin{equation}
\log_{\lambda}(1+t)=\sum_{n=1}^{\infty}\frac{\lambda^{n-1}(1)_{n,1/\lambda}}{n!}t^{n},\quad (\mathrm{see}\ [9]).\label{3}
\end{equation}
Note that 
\begin{displaymath}
	\lim_{\lambda\rightarrow 0}\log_{\lambda}(1+t)=\sum_{n=1}^{\infty}\frac{(-1)^{n-1}}{n}t^{n}=\log(1+t).
\end{displaymath}
It is well known that the Stirling numbers of the first kind are defined by 
\begin{equation}
(x)_{n}=\sum_{k=0}^{n}S_{1}(n,k)x^{k},\quad (n\ge 0),\quad (\mathrm{see}\ [1-24]),\label{4}
\end{equation}
where 
\begin{equation*}
(x)_{0}=1,\quad (x)_{n}=x(x-1)\cdots(x-n+1),\quad (n\ge 1).
\end{equation*}
The Stirling numbers of the second kind are given by 
\begin{equation}
x^{n}=\sum_{k=0}^{n}S_{2}(n,k)(x)_{k},\quad (n\ge 0),\quad (\mathrm{see}\ [1-24]). \label{5}
\end{equation}
Recently, the degenerate Stirling numbers of the first kind are defined by Kim-Kim as 
\begin{equation}
(x)_{n}=\sum_{k=0}^{n}S_{1,\lambda}(n,k)(x)_{k,\lambda},\quad (n\ge 0),\quad (\mathrm{see}\ [9]). \label{6}
\end{equation}
As the inversion formula of \eqref{6}, the degenerate Stirling numbers of the second kind are defined by 
\begin{equation}
(x)_{n,\lambda}=\sum_{k=0}^{n}S_{1,\lambda}(n,k)(x)_{k},\quad (\mathrm{see}\ [9]). \label{7}
\end{equation}
From \eqref{6} and \eqref{7}, we note that 
\begin{equation}
\frac{1}{k!}\Big(\log_{\lambda}(1+t)\Big)^{k}=\sum_{n=k}^{\infty}S_{1,\lambda}(n,k)\frac{t^{n}}{n!},\label{8}
\end{equation}
and 
\begin{equation*}
	\frac{1}{k!}\Big(e_{\lambda}(t)-1\Big)^{k}=\sum_{n=k}^{\infty}S_{2,\lambda}(n,k)\frac{t^{n}}{n!},\quad (\mathrm{see}\ [9,15]). 
\end{equation*}
For $r\ge 0$, the degenerate $r$-Stirling numbers of the first kind are defined by 
\begin{equation}
(x+r)_{n}=\sum_{k=0}^{n}S_{1,\lambda}^{(r)}(n+r,k+r)(x)_{k,\lambda},\quad (n\ge 0),\quad (\mathrm{see}\ [18,20,21]).\label{9}	
\end{equation}
In view of \eqref{5}, the degenerate $r$-Stirling numbers of the second kind are defined by 
\begin{equation}
(x+r)_{n,\lambda}=\sum_{k=0}^{n}{n+r \brace k+r}_{r,\lambda}(x)_{k},\quad (\mathrm{see}\ [18,21]).\label{10}
\end{equation}
The degenerate unsigned $r$-Stirling numbers of the first kind are given by 
\begin{equation}
\langle x+r\rangle_{n}=\sum_{k=0}^{n}{n+r \brack k+r}_{r,\lambda}\langle x\rangle_{k,\lambda},\quad (\mathrm{see}\ [18,20,21]),\label{11}	
\end{equation}
where the rising factorials and the degenerate rising factorials are respectively given by
\begin{align*}
&\langle x\rangle_{n}=x(x+1)\cdots(x+n-1),\quad (n\ge 1),\quad \langle x\rangle_{0}=1, \\
&\langle x\rangle_{n,\lambda}=x(x+\lambda)\cdots(x+(n-1)\lambda),\quad (n\ge 1),\quad \langle x\rangle_{0,\lambda}=1.	
\end{align*}
From \eqref{9}, \eqref{10} and \eqref{11}, we note that 

\begin{align}
&\frac{1}{k!}(1+t)^{r}\Big(\log_{\lambda}(1+t)\Big)^{k}=\sum_{n=k}^{\infty}S_{1,\lambda}^{(r)}(n+r,k+r)\frac{t^{n}}{n!}, \nonumber\\
&\frac{1}{k!}\Big(e_{\lambda}(t)-1\Big)^{k}e_{\lambda}^{r}(t)=\sum_{n=k}^{\infty}{n+r \brace k+r}_{r,\lambda}\frac{t^{n}}{n!},\label{12}\\
&\frac{1}{k!}\Big(-\log_{\lambda}(1-t)\Big)^{k}(1-t)^{-r}=\sum_{n=k}^{\infty}{n+r \brack k+r}_{r,\lambda}\frac{t^{n}}{n!},\quad (\mathrm{see}\ [18,20]). \nonumber
\end{align}
The degenerate Bell polynomials are defined by
\begin{equation}
e^{x(e_{\lambda}(t)-1)}=\sum_{n=0}^{\infty}\phi_{n,\lambda}(x)\frac{t^{n}}{n!},\quad (\mathrm{see}\ [15,16]). \label{13}	
\end{equation}
By \eqref{8} and \eqref{13}, we get 
\begin{equation}
\phi_{n,\lambda}(x)=\sum_{k=0}^{n}S_{2,\lambda}(n,k)x^{k},\quad (n\ge 0),\quad (\mathrm{see}\ [15,16]).\label{14}
\end{equation}
More generally, for $r\ge 0$, the degenerate $r$-Bell polynomials are given by 
\begin{equation}
e_{\lambda}^{r}(t)e^{x(e_{\lambda}(t)-1)}=\sum_{n=0}^{\infty}\phi_{n,\lambda}^{(r)}(x)\frac{t^{n}}{n!}. \label{15}	
\end{equation}
When $x=1$, $\phi_{n,\lambda}^{(r)}=\phi_{n,\lambda}^{(r)}(1)$ are called the degenerate $r$-Bell numbers. \par 
From \eqref{12} and \eqref{15}, we note that 
\begin{equation}
\phi_{n,\lambda}^{(r)}(x)=\sum_{k=0}^{n}{n+r \brace k+r}_{r,\lambda}x^{k},\quad (n\ge 0),\quad (\mathrm{see}\ [18,21]).\label{16}	
\end{equation}
The harmonic numbers are defined by 
\begin{equation}
H_{0}=1,\quad H_{n}=1+\frac{1}{2}+\cdots+\frac{1}{n},\quad (n\in\mathbb{N}). \label{17}	
\end{equation}
From \eqref{17}, we have 
\begin{equation}
-\frac{\log(1-x)}{1-x}=\sum_{n=1}^{\infty}H_{n}x^{n},\quad (\mathrm{see}\ [5,13,14,23]). \label{18}
\end{equation}
Recently, the degenerate harmonic numbers are introduced as 
\begin{equation}
H_{0,\lambda}=0	,\quad H_{n,\lambda}=\sum_{k=1}^{n}\frac{1}{\lambda}\binom{\lambda}{k}(-1)^{k-1},\quad (n\in\mathbb{N}),\quad (\mathrm{see}\ [13]). \label{19}
\end{equation}
Note that $\lim_{\lambda\rightarrow 0}H_{n,\lambda}=H_{n},\ (n\ge 1)$. \par 
From \eqref{19}, we have 
\begin{equation}
-\frac{1}{1-t}\log_{\lambda}(1-t)=\sum_{n=1}^{\infty}H_{n,\lambda}t^{n},\quad (\mathrm{see}\ [12,13]).\label{20}
\end{equation}
For $n\ge 0,\ r\ge 1$, the degenerate hyperharmonic numbers are given by 
\begin{equation}
H_{n,\lambda}^{(1)}=H_{n,\lambda},\quad H_{0,\lambda}^{(r)}=0,\quad(r\ge 1),\quad H_{n,\lambda}^{(r)}=\sum_{k=1}^{n}H_{n,\lambda}^{(r-1)},\quad (r\ge 2, n \ge 1),\quad (\mathrm{see}\ [14]).\label{21}	
\end{equation}
Thus, by \eqref{21}, we get 
\begin{equation}
-\frac{\log_{\lambda}(1-t)}{(1-t)^{r}}=\sum_{n=1}^{\infty}H_{n,\lambda}^{(r)}t^{n},\quad (\mathrm{see}\ [14]). \label{22}
\end{equation}
It is well known that the Fubini polynomials are defined by 
\begin{equation}
\frac{1}{1-x(e^{t}-1)}=\sum_{n=0}^{\infty}F_{n}(x)\frac{t^{n}}{n!},\quad (\mathrm{see}\ [5,6,8,17,19]).\label{23}	
\end{equation}
Note that $F_{n}(x)=\sum_{k=0}^{n}S_{2}(n,k)k!x^{k},\ (n\ge 0)$. \par 
The degenerate Fubini polynomials are defined by (see [17])
\begin{equation}
\frac{1}{1-x(e_{\lambda}(t)-1)}=\sum_{n=0}^{\infty}F_{n,\lambda}(x)\frac{t^{n}}{n!}.\label{24}
\end{equation}
Note that $\lim_{\lambda\rightarrow 0}F_{n,\lambda}(x)=F_{n}(x)$ and $F_{n,\lambda}
(x)=\sum_{k=0}^{n}S_{2,\lambda}(n,k)k!x^{k},\quad (n\ge 0)$. \par 

\vspace{0.1in}

\section{Some results on degenerate harmonic numbers and degenerate Fubini polynomials}
In this section, we assume that $r$ is a fixed nonnegative integer. \par 
For $f(x)=\sum_{n=0}^{\infty}a_{n}x^{n}\in\mathbb{C}[\![x]\!]$, we define 
\begin{equation}
f_{\lambda}(x)=\sum_{n=0}^{\infty}a_{n}(x)_{n,\lambda}\in \mathbb{C}[\![x]\!], \label{25}
\end{equation}
where $\lambda$ is any fixed real number.
First, we observe that 
\begin{align}
\bigg(x\frac{d}{dx}\bigg)_{m,\lambda}x^{r}f(x)&=\bigg(x\frac{d}{dx}\bigg)_{m,\lambda}\sum_{k=0}^{\infty}a_{k}x^{k+r}=\sum_{k=0}^{\infty}a_{k}(k+r)_{m,\lambda}x^{k+r}	\label{26}\\
 &=\sum_{k=0}^{\infty}a_{k}\sum_{l=0}^{m}{m+r \brace l+r}_{r,\lambda}(k)_{l}x^{k+r} \nonumber\\
 &=\sum_{l=0}^{m}{m+r \brace l+r}_{r,\lambda}\bigg(\sum_{k=0}^{\infty}a_{k}(k)_{l}x^{k-l}\bigg)x^{l+r}\nonumber \\
 &=\sum_{l=0}^{m}{m+r \brace l+r}_{r,\lambda}x^{l+r}f^{(l)}(x), \nonumber
\end{align}
where $f^{(l)}(x)=\big(\frac{d}{dx}\big)^{l}f(x)$. \par 
From \eqref{26}, we note that 
\begin{equation}
\begin{aligned}
\bigg(x\frac{d}{dx}\bigg)_{m-r,\lambda}x^{r}\bigg(\frac{d}{dx}\bigg)^{r}f(x)&=\sum_{l=0}^{m-r}{m \brace l+r}_{r,\lambda}x^{l+r}\bigg(\frac{d}{dx}\bigg)^{l+r}f(x) \\
&=\sum_{l=r}^{m}{m \brace l}_{r,\lambda}x^{l}\bigg(\frac{d}{dx}\bigg)^{l}f(x).
\end{aligned}\label{27}
\end{equation}
Therefore, by \eqref{26} and \eqref{27}, we obtain the following theorem. 
\begin{theorem}
	For $m,r\in\mathbb{Z}$ with $m\ge r\ge 0$, we have 
	\begin{equation*}
	\bigg(x\frac{d}{dx}\bigg)_{m,\lambda}x^{r}=\sum_{l=0}^{m}{m+r \brace l+r}_{r,\lambda}x^{l+r}\bigg(\frac{d}{dx}\bigg)^{l}. 	
	\end{equation*}
In particular, we also have
\begin{equation*}
	\bigg(x\frac{d}{dx}\bigg)_{m-r,\lambda}x^{r}\bigg(\frac{d}{dx}\bigg)^{r}=\sum_{l=r}^{m}{m \brace l}_{r,\lambda}x^{l}\bigg(\frac {d}{dx}\bigg)^{l}.
\end{equation*}
\end{theorem}
Assume that 
\begin{equation}
f(x)=\sum_{n=0}^{\infty}a_{n}x^{n}\in \mathbb{C}[\![x]\!],\quad g(x)=\sum_{n=0}^{\infty}b_{n}x^{n}\in \mathbb{C}[\![x]\!]. \label{28}
\end{equation}
By \eqref{28}, we easily get $a_{k}=\frac{f^{(k)}(0)}{k!},\ b_{k}=\frac{g^{(k)}(0)}{k!}$. \par 
From Theorem 1 and \eqref{28}, we have 
\begin{align}
&x^{r}\sum_{n=0}^{\infty}a_{n}\sum_{k=0}^{n}{n+r \brace k+r}_{r,\lambda}x^{k}g^{(k)}(x)=\sum_{n=0}^{\infty}a_{n}\bigg(x\frac{d}{dx}\bigg)_{n,\lambda}x^{r}g(x) \label{29}\\
&=\sum_{n=0}^{\infty}a_{n}\bigg(x\frac{d}{dx}\bigg)_{n,\lambda}\sum_{l=0}^{\infty}b_{l}x^{l+r}=\sum_{l=0}^{\infty}b_{l}\sum_{n=0}^{\infty}a_{n}\bigg(x\frac{d}{dx}\bigg)_{n,\lambda}x^{l+r} \nonumber \\
&=\sum_{l=0}^{\infty}b_{l}\sum_{n=0}^{\infty}a_{n}(l+r)_{n,\lambda}x^{l+r}=x^{r}\sum_{l=0}^{\infty}b_{l}\sum_{n=0}^{\infty}a_{n}(l+r)_{n,\lambda}x^{l}\nonumber\\
&=x^{r}\sum_{l=0}^{\infty}b_{l}f_{\lambda}(l+r)x^{l}.\nonumber
\end{align}
In particular, by Theorem 1 and \eqref{28}, we get 
\begin{align}
&\sum_{k=r}^{m}{m \brace k}_{r,\lambda}x^{k}\bigg(\frac{d}{dx}\bigg)^{k}g(x)=\bigg(x\frac{d}{dx}\bigg)_{m-r,\lambda}x^{r}\bigg(\frac{d}{dx}\bigg)^{r}\sum_{n=0}	^{\infty}\frac{g^{(n)}(0)}{n!}x^{n} \label{30} \\
&=\sum_{n=0}^{\infty}\frac{g^{(n)}(0)}{n!}\bigg(x\frac{d}{dx}\bigg)_{m-r,\lambda}x^{r}\bigg(\frac{d}{dx}\bigg)^{r}x^{n}=\sum_{n=0}^{\infty}\frac{g^{(n)}(0)}{n!}(n)_{r}\bigg(x\frac{d}{dx}\bigg)_{m-r,\lambda}x^{n}\nonumber\\
&=\sum_{n=0}^{\infty}\frac{g^{(n)}(0)}{n!}\binom{n}{r}r!(n)_{m-r,\lambda}x^{n}=\sum_{n=r}^{\infty}\frac{g^{(n)}(0)}{n!}\binom{n}{r}r!(n)_{m-r,\lambda}x^{n}.\nonumber
\end{align}
By \eqref{30}, we get 
\begin{equation}
\sum_{m=r}^{\infty}\frac{f^{(m)}(0)}{m!}\bigg(\sum_{k=r}^{m}{m \brace k}_{r,\lambda}x^{k}g^{(k)}(x)\bigg)=\sum_{n=r}^{\infty}\frac{g^{(n)}(0)}{n!}\binom{n}{r}r!\bigg(\sum_{m=r}^{\infty}\frac{f^{(m)}(0)}{m!}(n)_{m-r,\lambda}\bigg)x^{n}. \label{31}	
\end{equation}
Therefore, by \eqref{29} and \eqref{31}, we obtain the following theorem. 
\begin{theorem}
Let $f(x)=\sum_{n=0}^{\infty}a_{n}x^{n}$ and $g(x)=\sum_{k=0}^{\infty}b_{k}x^{k}$. Then we have 
\begin{displaymath}
\sum_{n=0}^{\infty}\frac{f^{(n)}(0)}{n!}\sum_{k=0}^{n}{n+r \brace k+r}_{r,\lambda}x^{k}g^{(k)}(x)=\sum_{n=0}^{\infty}\frac{g^{(n)}(0)}{n!}f_{\lambda}(n+r)x^{n}. 
\end{displaymath}	
In particular, we also have
\begin{displaymath}
\sum_{m=r}^{\infty}\frac{f^{(m)}(0)}{m!}\bigg(\sum_{k=r}^{m}{m \brace k}_{r,\lambda}x^{k}g^{(k)}(x)\bigg)=\sum_{n=r}^{\infty}\frac{g^{(n)}(0)}{n!}\binom{n}{r}r!\bigg(\sum_{m=r}^{\infty}\frac{f^{(m)}(0)}{m!}(n)_{m-r,\lambda}\bigg)x^{n}.
\end{displaymath}
\end{theorem}
From \eqref{16}, we recall that 
\begin{equation}
\phi_{n,\lambda}^{(r)}(x)=\sum_{k=0}^{n}{n+r \brace k+r}_{r,\lambda}x^{k},\quad (n\ge 0).\label{33}
\end{equation}
Let us take $g(x)=e^{x}$ and $f(x)=x^{m},\ (m\in\mathbb{N})$, Note that $f_{\lambda}(x)=(x)_{m,\lambda}$. \par 
From Theorem 2, we have 
\begin{equation}
e^{x}\sum_{k=0}^{m}{m+r \brace k+r}_{r,\lambda}x^{k}=\sum_{n=0}^{\infty}\frac{1}{n!}(n+r)_{m,\lambda}x^{n}.\label{34}	
\end{equation}
Thus, by \eqref{16} and \eqref{34}, we get 
\begin{equation}
\phi_{m,\lambda}^{(r)}(x)=\sum_{k=0}^{m}{m+r \brace k+r}_{r,\lambda}x^{k} 
=e^{-x}\sum_{n=0}^{\infty}\frac{1}{n!}(n+r)_{m,\lambda}x^{n}. \label{35} 
\end{equation}
In view of \eqref{15}, we consider the degenerate $r$-Fubini polynomials which are given by
\begin{equation}
\sum_{n=0}^{\infty}F_{n,\lambda}^{(r)}(x)\frac{t^{n}}{n!}=e_{\lambda}^{r}(t)\frac{1}{1-x(e_{\lambda}(t)-1)}.\label{36}	
\end{equation}
By \eqref{12} and \eqref{36}, we get 
\begin{align}
&\sum_{n=0}^{\infty}F_{n,\lambda}^{(r)}(x)\frac{t^{n}}{n!}=e_{\lambda}^{r}(t)\frac{1}{1-x(e_{\lambda}(t)-1)}\label{37}\\
&=\sum_{k=0}^{\infty}x^{k}\frac{k!}{k!}e_{\lambda}^{r}(t)\big(e_{\lambda}(t)-1\big)^{k} 
=\sum_{n=0}^{\infty}\bigg(\sum_{k=0}^{n}{n+r \brace k+r}_{r,\lambda}x^{k}k!\bigg)\frac{t^{n}}{n!}. \nonumber
\end{align}
Comparing the coefficients on both sides of \eqref{37},  we have 
\begin{equation}
F_{n,\lambda}^{(r)}(x)=\sum_{k=0}^{n}{n+r \brace k+r}_{r,\lambda}k!x^{k}, \quad (n\ge 0). \label{38}
\end{equation}
Let us take $g(x)=\frac{1}{1-x}$ in Theorem 2. Then we have 
\begin{equation}
\sum_{n=0}^{\infty}\frac{f^{(n)}(0)}{n!}\sum_{k=0}^{n}{n+r \brace k+r}_{r,\lambda}x^{k}\frac{k!}{(1-x)^{k+1}}=\sum_{n=0}^{\infty}f_{\lambda}(n+r)x^{n}.\label{39}
\end{equation}
Let us take $f(x)=x^{m}$ in \eqref{39}. Then we have 
\begin{equation}
\frac{1}{1-x}\sum_{k=0}^{m}{m+r \brace k+r}_{r,\lambda}k!\bigg(\frac{x}{1-x}\bigg)^{k}=\sum_{n=0}^{\infty}(n+r)_{m,\lambda}x^{n}. \label{40}	
\end{equation}
By \eqref{38} and \eqref{40}, we get 
\begin{equation}
\frac{1}{1-x}F_{m,\lambda}^{(r)}\bigg(\frac{x}{1-x}\bigg)=\sum_{n=0}^{\infty}(n+r)_{m,\lambda}x^{n}.\label{41}
\end{equation}
Therefore, by \eqref{41}, we obtain the following theorem. 
\begin{theorem}
For $m\ge 0$, we have 
\begin{displaymath}
	\frac{1}{1-x}F_{m,\lambda}^{(r)}\bigg(\frac{x}{1-x}\bigg)=\sum_{n=0}^{\infty}(n+r)_{m,\lambda}x^{n}.
\end{displaymath}
In particular, we also have
\begin{displaymath}
	F_{m,\lambda}^{(r)}=F_{m,\lambda}^{(r)}(1)=\sum_{n=0}^{\infty}(n+r)_{m,\lambda}\bigg(\frac{1}{2}\bigg)^{n+1}, 
\end{displaymath}
where $F_{m,\lambda}^{(r)}$ are called the $r$-Fubini numbers. 
\end{theorem}
By \eqref{7} and \eqref{8}, we easily get 
\begin{equation}
S_{2,\lambda}(n+1,k)=S_{2,\lambda}(n,k-1)+(k-n\lambda)S_{2,\lambda}(n,k), \label{42}
\end{equation}
where $n,k\in\mathbb{N}$ with $n\ge k$. \par 
From \eqref{24}, we note that 
\begin{align}
&x\frac{d}{dx}F_{n,\lambda}(x)+x\frac{d}{dx}\Big(xF_{n,\lambda}(x)\Big)-n\lambda F_{n,\lambda}(x)\label{43} \\
&=\sum_{k=1}^{n}kS_{2,\lambda}(n,k)k!x^{k}+x\sum_{k=0}^{n}S_{2,\lambda}(n,k)k!x^{k} \nonumber \\
&\quad +x^{2}\sum_{k=1}^{n}S_{2,\lambda}(n,k)k!kx^{k-1}-n\lambda\sum_{k=0}^{n}S_{2,\lambda}(n,k)k!x^{k}\nonumber\\
&=\sum_{k=1}^{n}kS_{2,\lambda}(n,k)k!x^{k}+\sum_{k=1}^{n+1}S_{2,\lambda}(n,k-1)(k-1)!x^{k} \nonumber \\
&\quad +\sum_{k=2}^{n+1}S_{2,\lambda}(n,k-1)(k-1)(k-1)!x^{k}-n\lambda\sum_{k=0}^{n}S_{2,\lambda}(n,k)k!x^{k} \nonumber \\
&=\sum_{k=1}^{n+1}\Big\{(k-n\lambda)S_{2,\lambda}(n,k)+S_{2,\lambda}(n,k-1)\Big\}k!x^{k} \nonumber \\
&=\sum_{k=1}^{n+1}S_{2,\lambda}(n+1,k)k!x^{k}=F_{n+1,\lambda}(x).\nonumber
\end{align}
Therefore, by \eqref{43}, we obtain the following theorem.
\begin{theorem}
	For $n\ge 0$, we have 
	\begin{equation*}
		x\frac{d}{dx}F_{n,\lambda}(x)+x\frac{d}{dx}\Big(xF_{n,\lambda}(x)\Big)-n\lambda F_{n,\lambda}(x)=F_{n+1,\lambda}(x).
	\end{equation*}
\end{theorem} 
From \eqref{3}, we note that 
\begin{align}
&\frac{1}{2!}\frac{d}{dx}\Big(\log_{\lambda}(1-x)\Big)^{2}\label{44} \\
&=-\frac{\log_{\lambda}(1-x)}{1-x}\Big(\lambda\log_{\lambda}(1-x)+1\Big)=-\frac{\lambda}{1-x}\Big(\log_{\lambda}(1-x)\Big)^{2}-\frac{\log_{\lambda}(1-x)}{1-x} \nonumber \\
&=-2!\lambda\sum_{n=2}^{\infty}{n+1 \brack 2+1}_{1,\lambda}\frac{x^{n}}{n!}+\sum_{n=1}^{\infty}H_{n,\lambda}x^{n} \nonumber \\
&=\sum_{n=1}^{\infty}\bigg(-2!\lambda {n+1 \brack 3}_{1,\lambda}+n!H_{n,\lambda}\bigg)\frac{x^{n}}{n!}. \nonumber
\end{align}
On the other hand, by \eqref{12}, we get 
\begin{align}
&\frac{1}{2!}\frac{d}{dx}\Big(\log_{\lambda}(1-x)\Big)^{2}=\frac{d}{dx}\sum_{n=2}^{\infty}S_{1,\lambda}(n,2)\frac{(-1)^{n}x^{n}}{n!} \label{45} \\
&=\sum_{n=2}^{\infty}S_{1,\lambda}(n,2)\frac{(-1)^{n}x^{n-1}}{(n-1)!}=\sum_{n=1}^{\infty}S_{1,\lambda}(n+1,2)(-1)^{n-1}\frac{x^{n}}{n!} \nonumber \\
&=\sum_{n=1}^{\infty}{n+1 \brack 2}_{\lambda}\frac{x^{n}}{n!},\nonumber	
\end{align}
where ${n \brack k}_{\lambda}=(-1)^{n-k}S_{1,\lambda}(n,k)$ are called the unsigned degenerate Stirling numbers of the first kind. \par 
By \eqref{12} and \eqref{22}, we get 
\begin{align}
-\frac{\log_{\lambda}(1-t)}{(1-t)^{r}} &=\sum_{n=1}^{\infty}H_{n,\lambda}^{(r)}t^{n} \label{46}\\
	&=\sum_{n=1}^{\infty}{n+r \brack r+1}_{r,\lambda}\frac{t^{n}}{n!}. \nonumber
\end{align}
From \eqref{44}, \eqref{45} and \eqref{46}, we have 
\begin{align}
&n!H_{n,\lambda}^{(r)}={n+r \brack r+1}_{r,\lambda},\label{47} \\
&n!H_{n,\lambda}=2!\lambda{n+1 \brack 3}_{1,\lambda}+{n+1 \brack 2}_{\lambda},\quad (n\ge 1).\nonumber
\end{align}
Therefore, by \eqref{47}, we obtain the following theorem. 
\begin{theorem}
For $n\in\mathbb{N}$, we have 
\begin{align*}
&n!H_{n,\lambda}^{(r)}={n+r \brack r+1}_{r,\lambda}, \\
&n!H_{n,\lambda}=2!\lambda{n+1 \brack 3}_{1,\lambda}+{n+1 \brack 2}_{\lambda},
\end{align*}
where ${n+1 \brack 2}_{\lambda}=(-1)^{n-1}S_{1,\lambda}(n+1,2)$ are the unsigned degenerate Stirling numbers of the first kind.
\end{theorem}
\begin{remark}
Note that Theorem 5 shows the following which holds between the degenerate 1-Stirling numbers of the first kind and the unsigned degenerate Stirling numbers of the first kind.
\begin{equation*}
{n+1 \brack 2}_{1,\lambda}=2!\lambda{n+1 \brack 3}_{1,\lambda}+{n+1 \brack 2}_{\lambda}.
\end{equation*}
\end{remark}
Let $g(t)=-\frac{1}{1-t}\log_{\lambda}(1-t)$. Then, by \eqref{19}, we get 
\begin{align}
g^{\prime}(t)&=-\frac{1}{(1-t)^{2}}\log_{\lambda}(1-t)+\frac{1}{(1-t)^{2}}\Big(\lambda\log_{\lambda}(1-t)+1\Big)\label{48} \\
&=-\frac{(1-\lambda)}{(1-t)^{2}}\log_{\lambda}(1-t)+\frac{1}{(1-t)^{2}} \nonumber \\
&=\frac{1!}{(1-t)^{2}}\bigg(-\binom{1-\lambda}{1}\log_{\lambda}(1-t)+H_{1,\lambda}\bigg),\nonumber
\end{align}
\begin{align*}
g^{\prime\prime}(t)&=\frac{d}{dt}g^{\prime}(t)=-\frac{2(1-\lambda)}{(1-t)^{3}}\log_{\lambda}(1-t)+\frac{(1-\lambda)(\lambda\log_{\lambda}(1-t)+1)}{(1-t)^{3}}+\frac{2!}{(1-t)^{3}} \\
&=-\frac{(1-\lambda)(2-\lambda)}{(1-t)^{3}}\log_{\lambda}(1-t)+\frac{2!}{(1-t)^{3}}\bigg(1+\frac{1-\lambda}{2}\bigg) \\
&=\frac{2!}{(1-t)^{3}}\bigg(-\binom{2-\lambda}{2}\log_{\lambda}(1-t)+\sum_{k=1}^{2}\frac{1}{\lambda}\binom{\lambda}{k}(-1)^{k-1}\bigg) \\
&=\frac{2!}{(1-t)^{3}}\bigg(-\binom{2-\lambda}{2}\log_{\lambda}(1-t)+H_{2,\lambda}\bigg),
\end{align*}
and 
\begin{align*}
	g^{(3)}(t)&=\frac{d}{dt}g^{\prime\prime}(t)=-\frac{3(1-\lambda)(2-\lambda)}{(1-t)^{4}}\log_{\lambda}(1-t)\\
	&\quad +\frac{(1-\lambda)(2-\lambda)(\lambda\log_{\lambda}(1-t)+1)}{(1-t)^{4}} +\frac{3!}{(1-t)^{4}}\bigg(\frac{1-\lambda}{2}+1\bigg) \\
	&=\frac{3!}{(1-t)^{4}}\bigg(-\binom{3-\lambda}{3}\log_{\lambda}(1-t)+\sum_{k=1}^{3}\frac{1}{\lambda}\binom{\lambda}{k}(-1)^{k-1}\bigg) \\
	&=\frac{3!}{(1-t)^{4}}\bigg(-\binom{3-\lambda}{3}\log_{\lambda}(1-t)+H_{3,\lambda}\bigg). 
\end{align*}
Continuing this process, we have 
\begin{equation}
g^{(k)}(t)=\frac{k!}{(1-t)^{k+1}}\bigg(-\binom{k-\lambda}{k}\log_{\lambda}(1-t)+H_{k,\lambda}\bigg),\label{49}
\end{equation}
and 
\begin{displaymath}
	g^{(k)}(0)=k!H_{k,\lambda},\quad (k\in\mathbb{N}).
\end{displaymath}
Therefore, by \eqref{49}, we obtain the following theorem. 
\begin{theorem}
	Let $g(t)=-\frac{1}{1-t}\log_{\lambda}(1-t)$. For $k\in\mathbb{N}$, we have 
	\begin{align*}
		g^{(k)}(t)&=\bigg(\frac{d}{dt}\bigg)^{k}g(t)
		=\frac{k!}{(1-t)^{k+1}}\bigg(-\binom{k-\lambda}{k}\log_{\lambda}(1-t)+H_{k,\lambda}\bigg),
	\end{align*}
	and 
	\begin{displaymath}
		g^{(k)}(0)=k!H_{k,\lambda}.
	\end{displaymath}
\end{theorem}
From \eqref{20}, \eqref{22} and Theorem 6, we note that 
\begin{align}
&\sum_{n=k}^{\infty}n(n-1)\cdots(n-k+1)H_{n,\lambda}t^{n-k}=\frac{d^{k}}{dt^{k}}\bigg(-\frac{\log_{\lambda}(1-t)}{1-t}\bigg) \label{50} \\
&=\frac{k!}{(1-t)^{k+1}}\bigg(H_{k,\lambda}-\binom{k-\lambda}{k}\log_{\lambda}(1-t)\bigg)=k!\bigg(\frac{H_{k,\lambda}}{(1-t)^{k+1}}-\binom{k-\lambda}{k}\frac{\log_{\lambda}(1-t)}{(1-t)^{k+1}}\bigg) \nonumber \\
&=k!\bigg\{\sum_{n=0}^{\infty}\binom{n+k}{n}H_{k,\lambda}t^{n}+\binom{k-\lambda}{k}\sum_{n=1}^{\infty}H_{n,\lambda}^{(k+1)}t^{n}\bigg\}. \nonumber
\end{align}
Thus, by \eqref{50}, we get 
\begin{equation}
\begin{aligned}
	&\sum_{n=0}^{\infty}\binom{n+k}{k}k!H_{n+k,\lambda}t^{n} \\
	&=\sum_{n=0}^{\infty}k!\binom{n+k}{k}H_{k,\lambda}t^{n}+k!\binom{k-\lambda}{k}\sum_{n=1}^{\infty}H_{n,\lambda}^{(k+1)}t^{n}. 
\end{aligned}	\label{51}
\end{equation}
By \eqref{51}, we get 
\begin{equation}
\begin{aligned}
&\sum_{n=1}^{\infty}\binom{n+k}{k}k!H_{n+k,\lambda}t^{n}\\
&\quad =\sum_{n=1}^{\infty}\bigg\{k!\binom{n+k}{k}H_{k,\lambda}+k!\binom{k-\lambda}{k}H_{n,\lambda}^{(k+1)}\bigg\}t^{n}.	
\end{aligned}\label{52}
\end{equation}
From \eqref{52}, we obtain the following theorem.
\begin{corollary}
For $n,k\in\mathbb{N}$, we have 
\begin{displaymath}
	\binom{k-\lambda}{k}H_{n,\lambda}^{(k+1)}=\binom{n+k}{k}\Big(H_{n+k,\lambda}-H_{k,\lambda}\Big).
\end{displaymath}
\end{corollary}
Let us take $g(x)=-\frac{1}{1-x}\log_{\lambda}(1-x)$ in Theorem 2. Then we have 
\begin{equation}
\begin{aligned}
	&\sum_{n=0}^{\infty}\frac{f^{(n)}(0)}{n!}\sum_{k=0}^{n}{n+r \brace k+r}_{r,\lambda}\frac{k!}{1-x}\bigg(\frac{x}{1-x}\bigg)^{k}\bigg(-\binom{k-\lambda}{k}\log_{\lambda}(1-x)+H_{k,\lambda}\bigg)\\
	&=\sum_{n=0}^{\infty}H_{n,\lambda}f_{\lambda}(n+r)x^{n}. 
\end{aligned}\label{53}
\end{equation}
For $f(x)=x^{m}$, by \eqref{53}, we get 
\begin{align}
&\sum_{n=0}^{\infty}H_{n,\lambda}(n+r)_{m,\lambda}x^{n} \label{54} \\
&=\frac{1}{1-x}\sum_{k=0}^{m}{m+r \brace k+r}_{r,\lambda}k!\bigg(\frac{x}{1-x}\bigg)^{k}H_{k,\lambda}-\frac{\log_{\lambda}(1-x)}{1-x}\sum_{k=0}^{m}{m+r\brace k+r}_{r,\lambda}k!\bigg(\frac{x}{1-x}\bigg)^{k}\binom{k-\lambda}{k} \nonumber \\
&=\frac{1}{1-x}\sum_{k=0}^{m}{m+r \brace k+r}_{r,\lambda}k!\bigg(\frac{x}{1-x}\bigg)^{k}\bigg(H_{k,\lambda}-\binom{k-\lambda}{k}\log_{\lambda}(1-x)\bigg).\nonumber
\end{align}
Therefore, by \eqref{54}, we obtain the following theorem. 
\begin{theorem}
	For $m\ge 0$, we have 
	\begin{align*}
		&\sum_{n=0}^{\infty}H_{n,\lambda}(n+r)_{m,\lambda}x^{n} \\
		&=\frac{1}{1-x}\sum_{k=0}^{m}{m+r \brace k+r}_{r,\lambda}k!\bigg(\frac{x}{1-x}\bigg)^{k}\bigg(H_{k,\lambda}-\binom{k-\lambda}{k}\log_{\lambda}(1-x)\bigg).
	\end{align*}
\end{theorem}
\section{conclusion}
A general identity valid for any two formal power series was derived, where the degenerate $r$-Stirling numbers of the second kind appear. Applying that result to appropriate formal power series, we obtained an identity on the degenerate $r$-Fubini polynomials. We deduced a differential equation satisfied by the degenerate Fubini polynomials with the help of a recurrence relation for the degenerate Stirling numbers of the second kind. The degenerate hyperharmonic numbers were expressed in terms of the degenerate unsigned $r$-Stirling numbers of the first kind. We derived a relation between the degenerate hyperharmonic numbers and the degenerate harmonic numbers.  By applying the general identity with suitable formal power series, we obtained an identity expressing a power series with coefficients given by the product of the degenerate harmonic numbers and the degenerate falling factorials in terms of the degenerate $r$-Stirlings of the second, the degenerate harmonic numbers and the degenerate logarithms. \par
It is one of our future research projects to continue to study various degenerate versions of some degenerate special numbers and polynomials and to find their applications in physics, science and engineering as well as in mathematics.


\begin{thebibliography}{9} 
\bibitem{1}
Abramowitz, M.; Stegun, I. A. \emph{Handbook of mathematical functions with formulas, graphs, and mathematical tables.} For sale by the Superintendent of Documents. National Bureau of Standards Applied Mathematics Series, No. 55 U. S. Government Printing Office, Washington, D.C., 1964 xiv+1046 pp.
\bibitem{2}
Broder, A. Z. \emph{The $r$-Stirling numbers.} Discrete Math. \textbf{49} (1984), no. 3, 241-259.
\bibitem{3}
Brualdi, R. A. \emph{Introductory combinatorics.} Fifth edition. Pearson Prentice Hall, Upper Saddle River, NJ, 2010. xii+605 pp. ISBN: 978-0-13-602040-0; 0-13-602040-2 
\bibitem{4}
Carlitz, L. \emph{Degenerate Stirling, Bernoulli and Eulerian numbers.} Utilitas Math. \textbf{15} (1979), 51-88.
\bibitem{5}
Comtet, L. \emph{Advanced combinatorics.} The art of finite and infinite expansions. Revised and enlarged edition. D. Reidel Publishing Co., Dordrecht, 1974. xi+343 pp. ISBN: 90-277-0441-4.
\bibitem{6}
Conway, J. H. ; Guy, R. K.{} \emph{The book of numbers.} (English summary) Copernicus, New York, 1996. x+310 pp. ISBN: 0-387-97993-X.
\bibitem{7}
Djordjevic, G. B.; Milovanovic, G. V. \emph{Special classes of polynomials.} University of Nis, Faculty of Technology, Leskovac, 2014. 
\bibitem{8}
Dougherty, M.; McCammond, J. \emph{Geometric combinatorics of polynomials I: The case of a single polynomial.} J. Algebra \textbf{607} (2022), 106-138.
\bibitem{9}
Kim, D. S.; Kim, T. \emph{A note on a new type of degenerate Bernoulli numbers.} Russ. J. Math. Phys. \textbf{27} (2020), no. 2, 227-235.
\bibitem{10}
 Kim, H. K. \emph{Fully degenerate Bell polynomials associated with degenerate Poisson random variables.} Open Math. \textbf{19} (2021), no. 1, 284-296.
\bibitem{11}
Kim, H. K. \emph{Combinatorial identities degenerate $r$-Dowling-Lah polynomials and numbers arising from degenerate umbral calculus.} Adv. Stud. Contemp. Math. (Kyungshang) \textbf{32} (2022), No. 3, 303-324.
\bibitem{12}
Kim, T.; Kim, D. S. \emph{Some relations of two type 2 polynomials and discrete harmonic numbers and polynomials.} Symmetry 2020, \textbf{12}, 905; doi:10.3390/sym12060905
\bibitem{13}
 Kim, T.; Kim, D. S. \emph{On some degenerate differential and degenerate difference operators.} Russ. J. Math. Phys. \textbf{29} (2022), no. 1, 37-46.
\bibitem{14}
Kim, T.; Kim, D. S. \emph{Some identities on degenerate hyperharmonic numbers.} arXiv:2205.10010  
\bibitem{15}
Kim, T.; Kim, D. S. \emph{Some identities on degenerate Bell polynomials and their related identities.} Proc. Jangjeon Math. Soc. \textbf{25} (2022), no. 1, 1-11.
\bibitem{16}
Kim, T.; Kim, D. S.; Dolgy, D. V. \emph{On partially degenerate Bell numbers and polynomials.} Proc. Jangjeon Math. Soc. \textbf{20} (2017), no. 3, 337-345.
\bibitem{17}
Kim, T.; Kim, D. S.; Kim, H. K.; Lee, H. \emph{Some properties on degenerate Fubini polynomials.} Appl. Math. Sci. Eng. \textbf{30} (2022), no. 1, 235–248.
\bibitem{18}
Kim, T.; Kim, D. S.; Kwon, J.; Lee, H. \emph{Some identities involving degenerate $r$-Stirling numbers.} Proc. Jangjeon Math. Soc. \textbf{25} (2022), no. 2, 245-252.
\bibitem{19}
Kim, T.; Kim, D. S.; Lee, H.; Kwon, J. \emph{On degenerate generalized Fubini polynomials.} AIMS Math. \textbf{7} (2022), no. 7, 12227-12240.
\bibitem{20}
Kim, T.; Kim, D. S.; Lee, H.; Park, J.-W. \emph{A note on degenerate $r$-Stirling numbers.} J. Inequal. Appl. 2020, Paper No. 225, 12 pp.
\bibitem{21}
Kim, T.; Yao, Y.; Kim, D. S.; Jang, G.-W. \emph{Degenerate $r$-Stirling numbers and $r$-Bell polynomials.} Russ. J. Math. Phys. \textbf{25} (2018), no. 1, 44-58.
\bibitem{22}
Pyo, Sung-Soo \emph{Some identities of degenerate Fubini polynomials arising from differential equations.} J. Nonlinear Sci. Appl. \textbf{11} (2018), no. 3, 383-393.
\bibitem{23}
Roman, S. \emph{The umbral calculus.} Pure and Applied Mathematics, 111. Academic Press, Inc. [Harcourt Brace Jovanovich, Publishers], New York, 1984. x+193 pp. ISBN: 0-12-594380-6.
\bibitem{24}
Washington, L. C. Elliptic curves. \emph{Number theory and cryptography.} Second edition. Discrete Mathematics and its Applications (Boca Raton). Chapman \& Hall/CRC, Boca Raton, FL, 2008. xviii+513 pp. ISBN: 978-1-4200-7146-7; 1-4200-7146-7.

\end{thebibliography}
\end{document}